\documentclass[]{article}
\usepackage{amsmath}

\title{Explicit expressions for real roots of a quartic equation}
\author{Nino Krvavica}
\date{%
	University of Rijeka, Faculty of Civil Engineering \\
	nino.krvavica@uniri.hr \\[2ex]
	\today
}

\begin{document}

\maketitle

This short article presents explicit expressions for roots of a quartic equation that has all four real roots. Although a general expression for quartic roots is available on Wikipedia \cite{wikipedia}, an optimized and slightly shorter expression for only real roots is presented here. A derivation of closed-form solutions for real roots of a quartic presented in the first part of this paper is taken from \cite{krvavica2018analytical}. It is repeated here for completeness and clarity.

Let us consider a general normalized 4th order polynomial equation (quartic)
\begin{equation}
x^4 + ax^3 + bx^2 + cx + d = 0.
\label{eq:quartic}
\end{equation}
To find the analytical solution to roots of Eq.~(\ref{eq:quartic}), first the cubic term $x^3$ is eliminated and the general polynomial is converted into a so-called \textit{depressed quartic} by a change of variables. Following Ferrari's method \cite{abramowitz1972}, a substitution $x = y - a/4$ is introduced, which gives a depressed polynomial
\begin{equation}
y^4 + py^2 + qy + r = 0,
\label{eq:depressed}
\end{equation}
where
\begin{equation}
p = b - 6\left(\frac{a}{4}\right)^2,
\end{equation}
\begin{equation}
q = c - 2b\left(\frac{a}{4}\right) + 8\left(\frac{a}{4}\right)^3,
\end{equation}
\begin{equation}
r = d - c \left(\frac{a}{4}\right) + b\left(\frac{a}{4}\right)^2 -3\left(\frac{a}{4}\right)^4. 
\end{equation}

The depressed polynomial can be rewritten as
\begin{equation}
\left(y^2 + \frac{p}{2}\right)^2 = -qy + \frac{p^2}{4} - r.
\label{eq:derpessed1}
\end{equation}
Next, expression $2zy^2 + zp + z^2$ is added to both sides of Eq.~(\ref{eq:derpessed1}), which after some regrouping gives
\begin{equation}
\left( y^2 + \frac{p}{2} + z \right)^2 = 2zy^2 -qy + z^2 + zp + \frac{p^2}{4} - r.
\label{eq:depressed2}
\end{equation}
When $z$ is chosen to be any non-zero root $z_0$ of the so-called \textit{resolvent cubic} equation
\begin{equation}
8z^3 + 8pz^2 + (2p^2 - 8r)z - q^2 = 0,
\label{eq:cubic1}
\end{equation}
the right-hand side of Eq.~(\ref{eq:depressed2}) can be written as a perfect square; therefore, Eq.~(\ref{eq:depressed2}) becomes
\begin{equation}
\left( y^2 + \frac{p}{2} + z_0 \right)^2 = \left(y\sqrt{2z_0} - \frac{q}{2\sqrt{2z_0}}\right)^2.
\label{eq:depressed3}
\end{equation}
And finally, Eq.~(\ref{eq:depressed3}) can be written as a factorized quadratic equation
\begin{equation}
\left( y^2 + \sqrt{2z_0}y + \frac{p}{2} + z_0 - \frac{q}{2\sqrt{2z_0}}\right)\left( y^2 - \sqrt{2z_0}y + \frac{p}{2} + z_0 + \frac{q}{2\sqrt{2z_0}}\right) = 0,
\end{equation}
which is easily solved by a quadratic formula.

Therefore, the solutions to the roots of the general quartic Eq.~(\ref{eq:quartic}) are given by
\begin{equation}
x_{1,2} = -\frac{a}{4} - \frac{1}{2}\sqrt{2z_0} \pm \frac{1}{2}\sqrt{- \left( 2p + 2z_0 - \frac{2q}{\sqrt{2z_0}} \right)},
\label{eq:solution1a}
\end{equation}
\begin{equation}
x_{3,4} = -\frac{a}{4} + \frac{1}{2}\sqrt{2z_0} \pm \frac{1}{2}\sqrt{- \left( 2p + 2z_0 + \frac{2q}{\sqrt{2z_0}} \right)}.
\label{eq:solution2b}
\end{equation}

For a general normalized 3rd order polynomial equation (cubic)
\begin{equation}
x^3 + \alpha x^2 + \beta x + \gamma = 0,
\label{eq:cubic}
\end{equation}
a real solution is given by Cardano's formula \cite{abramowitz1972}
\begin{equation}
x_0 = s_1 + s_2 - \frac{\alpha}{3},
\label{eq:x0}
\end{equation}
with
\begin{equation}
s_1 = \sqrt[3]{R + \sqrt{R^2 + Q^3 }},
\label{eq:S}
\end{equation}
\begin{equation}
s_2 = \sqrt[3]{R - \sqrt{R^2 + Q^3}},
\end{equation}
where
\begin{equation}
Q = \frac{3 \beta - \alpha^2}{9},
\end{equation}
\begin{equation}
R = \frac{9\alpha \beta - 27 \gamma - 2 \alpha^3}{54}.
\end{equation}

To eliminate redundant divisions and optimize computation of Eq.~(\ref{eq:solution1a}) and (\ref{eq:solution2b}), the root of the resolvent cubic equation is expressed via
\begin{equation}
2z_0 =  \frac{1}{3} \left( S + \frac{\Delta_0}{S}  - 2p \right),
\label{eq:cubic_Solution}
\end{equation}
where
\begin{equation}
S = 6s =  \sqrt[3]{\frac{\Delta_1 + \textrm{sign}(\Delta_1) \sqrt{\Delta_1^2 - 4\Delta_0^3 }}{2}},
\label{eq:S_sol}
\end{equation}
\begin{equation}
\Delta_0 = - 36 Q = b^2 + 12d - 3ac, 
\label{eq:Delta_0}
\end{equation}
\begin{equation}
\Delta_1 = 432 R = 27a^2d - 9abc + 2b^3 - 72bd + 27c^2.
\label{eq:Delta_1}
\end{equation}

Note that $\Delta_1^2 - 4 \Delta_0^3$, which is a much simpler expression for the discriminant of the resolvent cubic equation $\mathcal{D}_{cubic}$ and especially the discriminant of the quartic equation $\mathcal{D}_{quartic}$. Therefore, if $\Delta_1^2 - 4 \Delta_0^3 < 0$, three resolvent cubic roots are all real and the quartic roots are either all complex or all real. 
In this case, Eq.~(\ref{eq:cubic_Solution}) can be solved trigonometrically \cite{lambert1906}, which is computationally faster than computing the cube root required in Eq.~(\ref{eq:S_sol}):
\begin{equation}
z_0 =  \frac{1}{3} \left( \sqrt{\Delta_0} \cos \frac{\phi}{3}  - p \right),
\label{eq:cubic_trig}
\end{equation}
where
\begin{equation}
\phi = \arccos \left( \frac{\Delta_1}{2 \sqrt{\Delta_0^3}}\right).
\end{equation}

To summarize, a closed-form real solutions to the quartic equation can be simplified as follows:
\begin{equation}
\lambda_{1,4} = - \frac{a}{4} \mp \frac{ \sqrt{Z} + \sqrt{- A - Z \pm \frac{B}{\sqrt{Z}} } }{2} ,
\label{eq:solution1}
\end{equation}
\begin{equation}
\lambda_{2,3} = - \frac{a}{4} \mp \frac{ \sqrt{Z} - \sqrt{- A - Z \pm \frac{B}{\sqrt{Z}} } }{2} .
\label{eq:solution2}
\end{equation}
where
\begin{equation}
Z =  2z_0 = \frac{1}{3} \left( 2 \sqrt{\Delta_0} \cos \frac{\phi}{3}  - A \right),
\label{eq:Zcoeff}
\end{equation}
\begin{equation}
\phi = \arccos \left( \frac{\Delta_1}{2 \Delta_0 \sqrt{\Delta_0}}\right),
\label{eq:S_cubic}
\end{equation}
with
\begin{equation}
A = 2p = 2b - \frac{3a^2}{4},
\label{eq:Acoeff}
\end{equation}
\begin{equation}
B =2q = 2c - ab + \frac{a^3}{4} .
\label{eq:Bcoeff}
\end{equation}

Finally, by including Eqs.~(\ref{eq:Zcoeff})-(\ref{eq:Bcoeff}) and Eqs~(\ref{eq:Delta_0}) and (\ref{eq:Delta_1}) into Eqs.~(\ref{eq:solution1}) and (\ref{eq:solution2}), the fully explicit form of quartic roots is obtained:
\begin{equation}
\begin{aligned}
\lambda_{1,4} = & - \frac{a}{4} \mp \frac{1}{2} \sqrt{\frac{a^{2}}{4} - \frac{2 b}{3} + \frac{2}{3} \sqrt{- 3 a c + b^{2} + 12 d} \cos{\left (\frac{1}{3} \operatorname{acos}{\left (\frac{27 a^{2} d - 9 a b c + 2 b^{3} - 72 b d + 27 c^{2}}{\left(- 6 a c + 2 b^{2} + 24 d\right) \sqrt{- 3 a c + b^{2} + 12 d}} \right )} \right )}} \\
& \mp \frac{1}{2} \left[ \frac{a^{2}}{2} - \frac{4 b}{3} - \frac{2}{3} \sqrt{- 3 a c + b^{2} + 12 d} \cos{\left (\frac{1}{3} \operatorname{acos}{\left (\frac{27 a^{2} d - 9 a b c + 2 b^{3} - 72 b d + 27 c^{2}}{\left(- 6 a c + 2 b^{2} + 24 d\right) \sqrt{- 3 a c + b^{2} + 12 d}} \right )} \right)} \right. \\
& \pm \left. \frac{\frac{a^{3}}{4} - a b + 2 c}{\sqrt{\frac{a^{2}}{4} - \frac{2 b}{3} + \frac{2}{3} \sqrt{- 3 a c + b^{2} + 12 d} \cos{\left (\frac{1}{3} \operatorname{acos}{\left (\frac{27 a^{2} d - 9 a b c + 2 b^{3} - 72 b d + 27 c^{2}}{\left(- 6 a c + 2 b^{2} + 24 d\right) \sqrt{- 3 a c + b^{2} + 12 d}} \right )} \right )}}} \right]^{0.5}
\end{aligned}
\end{equation}

\begin{equation}
\begin{aligned}
\lambda_{2,3} = & - \frac{a}{4} \mp \frac{1}{2} \sqrt{\frac{a^{2}}{4} - \frac{2 b}{3} + \frac{2}{3} \sqrt{- 3 a c + b^{2} + 12 d} \cos{\left (\frac{1}{3} \operatorname{acos}{\left (\frac{27 a^{2} d - 9 a b c + 2 b^{3} - 72 b d + 27 c^{2}}{\left(- 6 a c + 2 b^{2} + 24 d\right) \sqrt{- 3 a c + b^{2} + 12 d}} \right )} \right )}} \\
& \pm \frac{1}{2} \left[ \frac{a^{2}}{2} - \frac{4 b}{3} - \frac{2}{3} \sqrt{- 3 a c + b^{2} + 12 d} \cos{\left (\frac{1}{3} \operatorname{acos}{\left (\frac{27 a^{2} d - 9 a b c + 2 b^{3} - 72 b d + 27 c^{2}}{\left(- 6 a c + 2 b^{2} + 24 d\right) \sqrt{- 3 a c + b^{2} + 12 d}} \right )} \right)} \right. \\
& \pm \left. \frac{\frac{a^{3}}{4} - a b + 2 c}{\sqrt{\frac{a^{2}}{4} - \frac{2 b}{3} + \frac{2}{3} \sqrt{- 3 a c + b^{2} + 12 d} \cos{\left (\frac{1}{3} \operatorname{acos}{\left (\frac{27 a^{2} d - 9 a b c + 2 b^{3} - 72 b d + 27 c^{2}}{\left(- 6 a c + 2 b^{2} + 24 d\right) \sqrt{- 3 a c + b^{2} + 12 d}} \right )} \right )}}} \right]^{0.5}
\end{aligned}
\end{equation}

\small
\bibliographystyle{abbrv} 
\bibliography{Qiqqa2BibTexExport}

\end{document}